\documentclass[12pt]{amsart}

\usepackage{amsmath,amssymb,latexsym,soul,cite,mathrsfs}

\usepackage{color,enumitem,graphicx}
\usepackage[colorlinks=true,urlcolor=blue,
citecolor=red,linkcolor=blue,linktocpage,pdfpagelabels,
bookmarksnumbered,bookmarksopen]{hyperref}
\usepackage[english]{babel}

\usepackage[left=2.9cm,right=2.9cm,top=2.8cm,bottom=2.8cm]{geometry}
\usepackage[hyperpageref]{backref}

\usepackage[colorinlistoftodos]{todonotes}
\makeatletter
\providecommand\@dotsep{5}
\def\listtodoname{List of Todos}
\def\listoftodos{\@starttoc{tdo}\listtodoname}
\makeatother

\numberwithin{equation}{section}

\begin{document}
	
	\title [On the origin of the Gibbons conjecture]{On the origin of the Gibbons conjecture}
	
	\author{Renan J. S. Isneri*}
	
	\pretolerance10000
	
	\address[Renan J. S. Isneri]
	{\newline\indent Unidade Acad\^emica de Matem\'atica
		\newline\indent 
		Universidade Federal de Campina Grande,
		\newline\indent
		58429-970, Campina Grande - PB - Brazil}
	\email{\href{renanisneri@mat.ufcg.edu.br}{renanisneri@mat.ufcg.edu.br}}
	
	\begin{abstract}
		\noindent  The term \emph{Gibbons conjecture} is widely used in connection with symmetry results for the Allen–Cahn equation. However, its origin is less transparent than its frequent citation suggests. In this note, we revisit its emergence, tracing it to a 1995 paper by Carbou and to subsequent developments in the literature. We argue that the attribution likely arose from informal communication rather than from a formally stated conjecture, illustrating how mathematical terminology may develop through transmission and collective usage.
	\end{abstract}

	
	
	
	\maketitle
	

\section{Introduction}

Many papers in nonlinear partial differential equations refer to the so-called \emph{Gibbons conjecture}. The statement appears frequently in the literature and is often used to describe a symmetry property of solutions to the Allen-Cahn equation. Yet a natural question arises: who actually formulated this conjecture?

To place this question in context, let us briefly recall a classical problem in the theory of nonlinear partial differential equations. 

The conjecture proposed by Ennio De Giorgi in the late 1970s has become one of the central symmetry problems in the study of entire solutions $u:\mathbb{R}^n \to \mathbb{R}$ of the Allen-Cahn equation
$$
-\Delta u + u^3 - u = 0 \quad \text{in}\quad \mathbb{R}^n,
$$
where $\Delta=\sum_{i=1}^{n}\partial_{x_i}^2$ denotes the Laplacian. This equation, introduced in the pioneering work of Allen and Cahn \cite{Allen}, admits the constant solutions $u \equiv -1$ and $u \equiv 1$, which play a fundamental role in its analysis. It arises naturally in the mathematical description of phase transitions, where the function $u$ represents different states of a material.

In its classical form, De Giorgi’s conjecture concerns bounded solutions that are monotone in one direction, for instance
$$
\partial_{x_n} u(x) > 0 \quad \text{for all } x \in \mathbb{R}^n.
$$
It predicts that such solutions must be one-dimensional, at least in low dimensions \cite{Giorgi}, meaning that $u$ depends only on a single variable: there exist a function $g:\mathbb{R}\to\mathbb{R}$, $a\in\mathbb{R}$ and a direction $\xi \in \mathbb{R}^{n}$ such that $u(x) = g(x \cdot \xi+a)$, where $x \cdot \xi$ denotes the usual inner product of $x$
and $\xi$. In other words, the level sets of $u$ are hyperplanes, and the solution varies only along one direction.

Over the past decades, De Giorgi’s conjecture has attracted considerable attention and has led to a rich body of results. It has revealed deep connections with minimal surface theory, variational methods, and geometric analysis. While substantial progress has been achieved in several dimensions (see, for example, \cite{Gui,Ambrosio,Pino,Gui1} and the references therein), the conjecture remains open in its full generality.

Beyond its intrinsic mathematical interest, the Allen-Cahn equation and the problems surrounding De Giorgi’s conjecture also arise in models from mathematical physics, including phenomena related to phase separation and cosmological domain walls (see, for instance, \cite{Gibbons}). 

Within this broader context, a closely related symmetry statement has frequently appeared in the literature under the name of the Gibbons conjecture.

Roughly speaking, instead of assuming monotonicity in one direction as in De Giorgi’s conjecture, one prescribes the behavior of the solution at infinity. More precisely, one considers bounded solutions $u:\mathbb{R}^n \to \mathbb{R}$ such that
$$
\lim_{x_n\to\pm\infty} u(x',x_n) = \pm 1 \quad\text{uniformly with respect to } x' \in\mathbb{R}^{n-1}.
$$
This condition means that, as one moves far in the $x_n$-direction, the solution becomes close to the constant solutions $-1$ and $1$, independently of the other variables. The conjecture asserts that any such solution must be one-dimensional and, in fact, is given (up to translation) by the explicit profile
$$
u(x) = \tanh\!\left(\frac{x_n-\alpha}{\sqrt{2}}\right),
$$
for some $\alpha \in \mathbb{R}$. In other words, the solution depends only on one variable and coincides with the classical transition between the two constant values. This statement is often viewed as a variant of De Giorgi’s conjecture, where the monotonicity assumption is replaced by an asymptotic condition at infinity.

In the classical semilinear setting, this conjecture has in fact been confirmed in all dimensions by different methods (see, for instance, \cite{Farina,Barlow,berestycki}). Nevertheless, the terminology has persisted and continues to appear in a variety of contexts, where it is used to describe symmetry and rigidity phenomena under asymptotic boundary conditions. In recent years, similar questions have been investigated for more general classes of equations, including quasilinear, nonlocal, and fully nonlinear models (see, for example, \cite{Esposito,FarinaValdinoci,BVDV}). These developments show that the ideas behind this conjecture remain active and continue to shape current research in nonlinear partial differential equations.

Despite its widespread use, the historical origin of this attribution is far from clear. The conjecture is commonly associated with Gary W. Gibbons, a theoretical physicist whose work influenced the study of equations of Ginzburg-Landau type, including the Allen-Cahn equation. However, no explicit formulation of this symmetry statement appears in his published work. This naturally leads to the question: how did this designation come to bear his name? 

In this note, we examine the emergence of the term Gibbons conjecture in the mathematical literature and clarify the circumstances under which this attribution arose. Our approach is based on tracing references and examining how the terminology appears and evolves across different sources. 

This case belongs to what may be regarded as the recent history of mathematics, where terminology and attribution often develop through informal exchanges, lectures, and the gradual consolidation of language within a research community. 

It is therefore somewhat surprising that a statement so frequently cited in the literature has no clearly identifiable original source.

We hope that this account will be of interest not only to specialists in nonlinear partial differential equations, but also to a broader audience interested in how mathematical ideas and terminology evolve.


\section{The Emergence of the Name}

Unlike De Giorgi’s conjecture, whose origin is well documented, the emergence of the so-called \emph{Gibbons conjecture} is much less clear. 

A careful examination of the literature suggests that the expression Gibbons conjecture can be traced back, at least in written form, to a 1995 paper by G. Carbou \cite{Carbou}, entitled \emph{Unicit\'e et minimalit\'e des solutions d’une \'equation de Ginzburg-Landau}. 

In that work, Carbou studies a particular class of solutions of the Allen-Cahn equation, namely minimizing solutions with prescribed behavior at infinity. Roughly speaking, these are solutions that connect the constant solutions $-1$ and $1$ as one moves far in a given direction. In this context, he emphasizes the role of the one-dimensional profile given by the hyperbolic tangent, which describes this transition.

In the course of this discussion, Carbou attributes a related uniqueness problem to Gary W. Gibbons, citing a private communication entitled {\it Topological Defects in Cosmology}. Notably, this communication does not appear as a published work.

According to Carbou, Gibbons had raised questions about the uniqueness, up to translations, of entire solutions in $\mathbb{R}^n$ with prescribed limits at infinity. These questions were motivated by models from cosmology involving equations of Ginzburg-Landau type.

It is worth noting, however, that Carbou does not explicitly use the expression Gibbons conjecture in his paper. Rather, he refers to Gibbons’ communication as a source of motivation for the uniqueness problem under consideration. The terminology itself appears to have been adopted subsequently by other authors.

In the years following Carbou’s work, the expression Gibbons conjecture began to appear explicitly in the literature. To the best of our knowledge, some of the earliest occurrences can be found in the works of Farina (1999) \cite{Farina} and Barlow, Bass, and Gui (2000) \cite{Barlow}, where the term is used in connection with symmetry results for solutions with prescribed behavior at infinity. In these works, the conjecture is established independently by different methods, contributing to the consolidation of the terminology in the literature. In all these cases, the reference ultimately traces back, either directly or indirectly, to Carbou’s 1995 paper.

The terminology has also remained in use in more recent works, including \cite{Esposito} (2022), \cite{GMZ} (2024), and \cite{CM} (2025), where the expression Gibbons conjecture even appears in the titles of the papers. This suggests that it has become firmly established in the language of the subject. 

Thus, although Carbou did not explicitly formulate the conjecture under this name, his paper appears to have played a central role in the emergence and dissemination of the expression, which has since become standard in the literature despite the absence of a clearly identifiable original source.

The reference to Gibbons cited by Carbou contains no bibliographical details (such as year, venue, or publisher) and does not correspond to any identifiable published manuscript. Moreover, a review of Gibbons’ published works reveals no article or book in which the symmetry statement now commonly referred to as the Gibbons conjecture is explicitly formulated. This suggests that the attribution likely originates from informal communication rather than from a formally documented source.

An illuminating remark appears in a 2009 paper by Gibbons and collaborators entitled \emph{The Bernstein Conjecture, Minimal Cones, and Critical Dimensions} \cite{Gibbons}. In a footnote (p.~5), he observes:
\begin{quote}
	\it ``This is known, for reasons that are only partially clear to GWG (one of the present authors), as the Gibbons Conjecture and has been proved by a number of people, but not by GWG.''
\end{quote}
A similar comment appears in a lecture delivered by Gibbons in Florence in May 2009, entitled \emph{The Bernstein Conjecture and the Fate of the 8-Brane}, where he again refers to the attribution in comparable terms (see footnote on p. 28 of the lecture notes). 

These remarks suggest that even Gibbons himself did not clearly identify the circumstances under which the conjecture came to bear his name, and indicate a certain degree of uncertainty regarding the origin of this attribution.

It is worth noting that, in the 2009 paper mentioned above, the Allen-Cahn equation is discussed in connection with problems from mathematical physics, including domain walls and minimal surfaces. This broader context helps explain why Carbou may have referred to Gibbons’ work as a source of motivation for the uniqueness problem he studied. 

Although the specific formulation now known as the \emph{Gibbons conjecture} cannot be traced to a written statement by Gibbons, his discussions of related physical models indicate that he was indeed engaged with problems closely connected to the symmetry issues later emphasized in the Partial Differential Equations literature. In this sense, the attribution reflects not a formal authorship of the conjecture as currently stated, but rather a genuine intellectual influence within a shared line of inquiry.

Additional clarification was obtained through correspondence with both Carbou and Gibbons. Carbou recalled that the private communication attributed to G. W. Gibbons originated from a seminar delivered by Gibbons in Paris, which he and his doctoral advisor attended and which motivated the formulation of the uniqueness problem in his 1995 article. 

Gibbons, for his part, indicated that he has no recollection of having written a manuscript entitled \emph{Topological Defects in Cosmology} and is unaware of any such document in his files. He suggested that the reference most likely corresponds to an unpublished lecture and also remarked that he himself is uncertain why the conjecture came to bear his name.

Taken together, these accounts strongly support the view that the problem originated in informal discussions and lectures, rather than in a formally stated conjecture in the literature.

What emerges from this reconstruction is not a case of mistaken attribution, but rather a gradual process through which a problem discussed in lectures and shaped by physical motivations became linked to a particular name. Over time, through repeated citation and collective usage, this association became established as standard terminology within the mathematical community.


\section{Attribution and Mathematical Memory}

The evidence discussed in this note does not point to a single written source in which Gibbons explicitly formulated the conjecture that now bears his name. At the same time, it also suggests that the attribution did not arise arbitrarily. Carbou’s 1995 paper, later remarks by Gibbons, and the surrounding context all indicate a genuine scientific connection, likely rooted in lectures, seminars, and shared interests in problems motivated by mathematical physics.

It is also noteworthy that in subsequent publications Gibbons cites Carbou’s 1995 paper, see for instance in \cite{Gibbons1}, and discusses related rigidity and Bernstein-type problems. This suggests an ongoing awareness of the analytical developments surrounding the conjecture and reinforces the existence of a scientific dialogue, even if indirect, between the physical motivations and the later formulation of the conjecture.

From this perspective, the name Gibbons conjecture appears not as the result of a single formal act of authorship, but rather as the outcome of a gradual process of transmission. A problem first discussed in one context and later reformulated in another eventually became established through repeated scholarly usage. In this sense, attribution in mathematics may reflect intellectual influence as much as formal authorship.

This episode illustrates a broader feature of mathematical culture. Names associated with conjectures and theorems do not always originate in explicit publication; they may instead emerge through lectures, informal exchanges, and subsequent consolidation in the literature. Such attributions are often natural, but they can also obscure the more complex paths through which ideas develop and gain acceptance.

Revisiting cases of this kind helps clarify not only a specific attribution, but also the mechanisms by which mathematical terminology becomes stabilized within the community. 

The history of the Gibbons conjecture illustrates how mathematical terminology may emerge not from a single formal statement, but from a combination of informal communication, subsequent citation, and collective usage. More broadly, it shows how names in mathematics may arise from influence and circulation as much as from formal statement.

%
%
%
%
%
%
%


\end{document}